\documentclass[lettersize,journal]{IEEEtran}
\usepackage{amsmath,amsfonts}
\usepackage{algorithmic}
\usepackage{algorithm}
\usepackage{array}
\usepackage[caption=false,font=normalsize,labelfont=sf,textfont=sf]{subfig}
\usepackage{textcomp}
\usepackage{stfloats}
\usepackage{url}
\usepackage{verbatim}
\usepackage{graphicx}
\usepackage{cite}
\usepackage{color,amsthm}
\usepackage{nomencl}
\usepackage{graphicx}
\usepackage[version=4]{mhchem}
\usepackage{siunitx}
\usepackage{longtable,tabularx}
\usepackage{comment}
\usepackage{caption}
\usepackage{multicol}
\usepackage{textcomp}
\usepackage{float}
\usepackage{wrapfig}
\usepackage{varwidth}
\usepackage{booktabs}
\usepackage{makecell}
\usepackage{nccmath}
\usepackage{lipsum}
\usepackage{mathtools}
\usepackage{cuted}
\usepackage[inline]{enumitem}
\usepackage{amssymb}

\begin{document}

\title{To charge in-flight or not: an inquiry into parallel-hybrid 
electric aircraft configurations\\ via optimal control}

\author{Mengyuan Wang,~\IEEEmembership{Student Member, IEEE} and Mehran Mesbahi,~\IEEEmembership{Fellow, IEEE}
\thanks{Mengyuan Wang and Mehran Mesbahi are with the William E. Boeing Department of Aeronautics and Astronautics, University of Washington, Seattle, USA,
        {\tt\small mywang10@uw.edu, mesbahi@uw.edu.}}
}



\maketitle

\begin{abstract}
We examine two configurations for parallel hybrid electric aircraft, one with, and one without, a mechanical connection between the engines and the electric motors. For this two designs, we then review the power allocation problem in the context of aircraft energy management for a 19-seat conceptual Hybrid Electric Aircraft. We then represent the original optimal control problem as a finite-dimensional optimization and validate the second-order sufficient conditions for global optimality of the obtained solution. This is then followed by a sensitivity analysis of the fuel consumption on the initial aircraft weight and flight endurance. Our simulation and theoretical results clarify the limited benefit of charging the battery in-flight for this class of hybrid electric aircraft to reduce $\mbox{CO}_2$ emissions.
\end{abstract}


\section{Introduction}\label{Sec_Introduction}

Aircraft electrification highly depends on the development of more efficient battery technologies, as the current energy density of most advanced batteries falls short of the desired levels for commercial aircraft. Consequently, developing Hybrid-Electric Aircraft (HEA) has become a viable strategy for next-generation aviation to satisfy reduced $\mbox{CO}_2$ emissions requirements in the coming decades. 

The existing literature on hybrid architectures for aircraft propulsion~\cite{ye2021review,wall2017survey,friedrich2015hybrid,brelje2019electric} mostly borrows principles from hybrid electric ground vehicles, often assuming that the gas turbine and electric motor should be connected mechanically, allowing the turbine to charge the battery during flight. This connective mechanism has several potential benefits: (a) The engine can run in its most efficient region to conserve fuel by using the electric motor as an extra load, and (b) The electric motor can work in``one-engine-inoperative'' (OEI) mode by charging the battery to its fully charged state after the climbing phase. These benefits are compatible with series, parallel, and series-parallel configurations. On the other hand, MagniX--a manufacturer of electric propulsion systems for electric aircraft--has proposed a novel parallel hybrid structure in 2021, in which the engine and electrical paths are independent (without a mechanical connection in place). In this structure, the electric motor and the engine work together during the climbing phase, and the engine handles the cruise phase and the descent phase. In this configuration, one can bypass the use of the engine and transition the HEA to an All-Electric Aircraft (AEA)--when the energy density of the battery reaches the desired level. In this paper, these two structures are compared in the context of the aforementioned benefits, and specifically, through the lens of optimal power allocation.

With two available energy sources for the propulsion system, the importance of power allocation algorithms cannot be overstated for HEA. In this directon, Leite and Voskuijl~\cite{leite2020optimal} utilized dynamic programming to obtain an optimal energy management scheme for HEA for a given flight profile; Doff-Sotta~{\em et al.}~\cite{doff2020optimal} on the other hand has developed a convex formulation for optimal HEA energy management. The latter work, in particular, showed that optimal control strategy significantly reduces fuel consumption as compared with heuristic methods. Donateo~{\em et al.}~\cite{donateo2018applying} in the meantime, applied dynamic programming for energy management for a lightweight rotorcraft, where the corresponding optimal control was formalized in terms of a shortest path problem on a graph, followed by the application of Dijkstra algorithm for its solution.

Although there are a few manufactured HEA on the market~\cite{martini2011world, wittmann2013flying,langelaan2013green,caujolle2017airbus}, publicly accessible technical data and characteristics for these aircraft are rather limited. In the meantime, there are a number of studies that have explored the design space of HEA. For example, De Vries~{\em et al.}~\cite{de2019preliminary} provided a comprehensive approach for the conceptual design process of HEA. Multiple hybrid-electric structures, including the all-electric propulsion system, have also been considered. Isikveren {\em et al.} examined the design of vertical take-off and landing aircraft~\cite{isikveren2014pre}. Finger {\em et al.} compared two design approaches utilizing an existing 19-seat conventional aircraft as the baseline model~\cite{finger2020comparison}.


\IEEEpubidadjcol 
In this work, we first review the power allocation problem for a given flight profile, formalized in terms of minimizing the corresponding fuel consumption by solving an optimal control problem (OCP). Subsequently, using an approximate linear fuel flow rate model, this OCP is represented as a finite-dimensional optimization problem. This formulation of the power allocation problem is then applied to a 19-seat conceptual aircraft model, where the second-order sufficient conditions are also examined for the corresponding optimization problem. In this case, the optimal power allocation strategy is tailored specifically for this mid-size aircraft with turboprop engines. Lastly, two hybrid structures (with and without mechanical connection in the propulsion system) are compared in the context of fuel optimized power allocation, allowing us to evaluate the utility of the mechanical connection to charge the battery in-flight.

The remainder of the paper is organized as follows. In Section~\ref{Sec_AircraftModel}, two parallel hybrid-electric architectures and the technical data of the conceptual aircraft model are introduced. In Section~\ref{Sec_FlightDynamics}, the longitudinal aircraft dynamics and simplifications made for the formulation of optimal control problems are provided; the power relation in the propulsion system is also discussed. In Section~\ref{Sec_PowerAllocation}, the previously proposed optimal control problem is reviewed and re-written based on a linear approximation of the engine's flow rate. Subsequently, this optimal control setup is formulated as a finite-dimensional optimization problem. Using the power allocation optimization formulation in Section~\ref{Sec:comparison}, and the corresponding sensitivity analysis, the fuel consumption between charging and not charging the battery during the cruise phase is then compared. This is then followed by evaluating the aircraft fuel consumption with different battery configurations (leading to different take-off weights). These two facets shed light on our evaluation for the question posed in the title of the manuscript. Conclusions and future work are then provided in Section~\ref{Sec:conclusion}.

\section{Hybrid-electric architectures and a conceptual aircraft model}\label{Sec_AircraftModel}
In this section, two parallel hybrid-electric architectures are discussed in Section~\ref{Subsec_HybridStructures}; technical data for the conceptual aircraft model and components in the propulsion system are provided in Section~\ref{Subsec_Aircraft model}.
\subsection{Hybrid-electric architectures}\label{Subsec_HybridStructures}
For a conventional aircraft, relying on engines for propulsion, it is challenging to devise means of operating the engine in its peak efficiency regime; this is primary due to  considerable differences in power requirements for the climb and the cruise phases of a typical flight profile. One of the primary advantages of hybrid electric aircraft is the ability to downsize the engine and use the electric motor as a supplemental propulsion, allowing the engine to operate at its most efficient regime during flight. This point of view when adopted for developing aircraft hybrid-electric architectures, is often combined with design principles from ground vehicles; as such, HEA designers often assume a mechanical connection in the propulsion system. The most commonly used architectures include series, parallel, and series-parallel configurations. The parallel architecture is the most popular configuration due to its relatively higher efficiency and lower weight (in comparison with the series configuration), as well as its simpler energy management (as compared with the series-parallel configuration). In general, a mechanical connection in the propulsion system enables conventional engines to work more efficiently, either by decoupling the engine from the propeller (in the series configuration) or by taking the electric motor as an extra load for the engine (in the parallel configuration). Nevertheless, to the best of our knowledge, this potential advantage of having a mechanical connection in HEA has not been validated, neither theoretically nor experimentally. 

The mechanical connection between the conventional engine and the electric motor often
necessitates a complicated clutch/gearing mechanism and advanced control algorithm for the hybrid propulsion system. In addition, charging the battery during flight introduces new challenges, such as the need for more sophisticated cooling systems and the possibility of battery pack degradation. 

MagniX has recently proposed a novel hybrid structure for aircraft electric propulsion. In particular, the Magnix design has the conventional engine and electric paths {\em mechanically independent}, with completely separate control systems for two propulsion paths. With this architecture, the hybrid electric propulsion system can still be advantageous as it allows efficient engine operation. In addition, the absence of mechanical connections reduces the aircraft's total weight and may improve safety due to less complexity in the mechanical system.

\begin{figure}[h]
    \centering
    \includegraphics[width=0.45\textwidth]{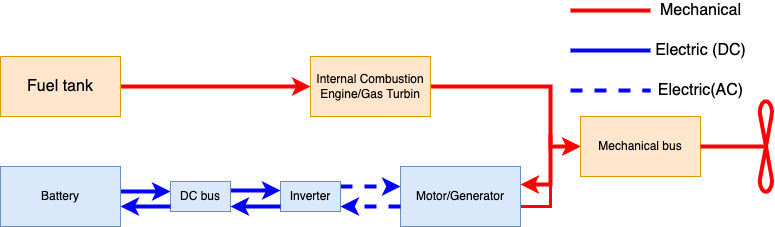}
    \caption{Connected parallel hybrid architecture}
    \label{Fig_ConnectedParallel}
\end{figure}

In this paper, we examine the standard parallel architecture (referred to as the connected configuration), and the MagniX architecture, which is called the independent configuration. In the connected configuration, the internal combustion engine/gas turbine and the electric motor are both connected to the propeller, as illustrated in Fig.~\ref{Fig_ConnectedParallel}. As a result, this configuration supports the following operating modes: 
\begin {enumerate*} [label=\itshape\alph*\upshape)]
\item engine-alone mode, \item motor-alone mode, \item combined mode (in which the engine and the motor drive the propeller simultaneously) and \item power split mode (in which the engine drive the propeller and charge the battery). 
\end {enumerate*}

\begin{figure}[h]
    \centering
    \includegraphics[width=0.45\textwidth]{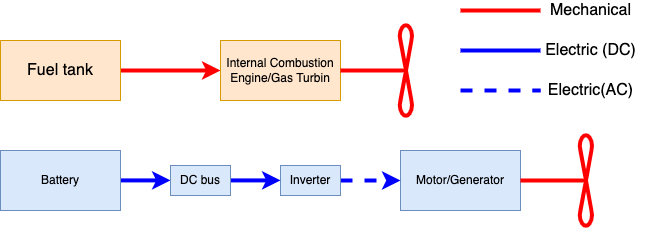}
    \caption{Independent parallel hybrid architecture}
    \label{Fig_IndependentParallel}
\end{figure}

In an independent configuration in Fig.~\ref{Fig_IndependentParallel}, the propulsion system retains the first three operating modes; the only difference is that the engine cannot charge the battery during flight. This paper provides a complete comparison between the fuel consumption of these two architectures using the power allocation problem formulation and the corresponding algorithm and sensitivity analysis~\cite{wang2021power}. 

Before delving into the power allocation problem formulation and its consequences for the main theme of this manuscript, we provide an overview of the aircraft model and its key components and subsystems of relevance to the power allocation problem and subsequent analysis.
\subsection{HEA propulsion system model}\label{Subsec_Aircraft model} 
We use the conceptual HEA propulsion system model developed by Finger, {\em et al.}~\cite{finger2020comparison} for subsequent analysis. Finger, {\em et al.} considered an existing 19-seat commuter aircraft as the baseline and examine two approaches for the HEA preliminary design. The conventional aircraft was re-sized to accommodate the hybrid electric propulsion technology. Some of the key characteristics of one of the configurations used in~\cite{finger2020comparison} are shown in Table~\ref{Table: HEA data}. Note that parameters with $^*$ are only estimated values.
\begin{table}[h]
    \centering
    \caption{Conceptual HEA technical data, including estimated parameters}
	\begin{tabular}{ll}
			\toprule
			Parameters  & Values\\
			\midrule
        Wing area & 32 $\mbox{m}^2$\\
        Maximum take-off weight & 6385 $\,\mbox{kg}$\\
        Cruise altitude & $3000\, \mbox{m}$\\
        Battery capacity$^*$ & 400$\,\mbox{Ah}$\\
        Battery pack nominal voltage$^*$ & 540$ \, \mbox{V}$\\
        Propeller efficiency$^*$ & 0.7\\
        Zero-lift drag coefficient$^*$ & 0.024\\
        Induced drag coefficient$^*$ & 0.056\\
        \bottomrule
		\end{tabular}
    \label{Table: HEA data}
\end{table}

In~\cite{finger2020comparison}, the authors did not specify the engine model, battery model, or electric motor model; instead, the reference weight and some characteristics for these components were provided. We have chosen the following models for the propulsion system--with some parameters slightly different from those in~\cite{finger2020comparison}:

\begin{enumerate}
    \item \textbf{Turboprop engine}: the PT6A-135A engine by Pratt \& Whitney is selected for our analysis; the weight of this engine is $153\mbox{kg}$; its maximum output power is $560\mbox{kW}$; and its average specific fuel consumption (SFC) is $356\mbox{g/kW/h}$. In this paper, it is assumed that the engine runs in a fixed rotational speed; the SFC and the fuel flow rate can be represented as functions of the output power of the engine, shown in Fig.~\ref{Fig_Preliminary_fuel}.
    \begin{figure}
        \subfloat[]{\includegraphics[width=1.7in]{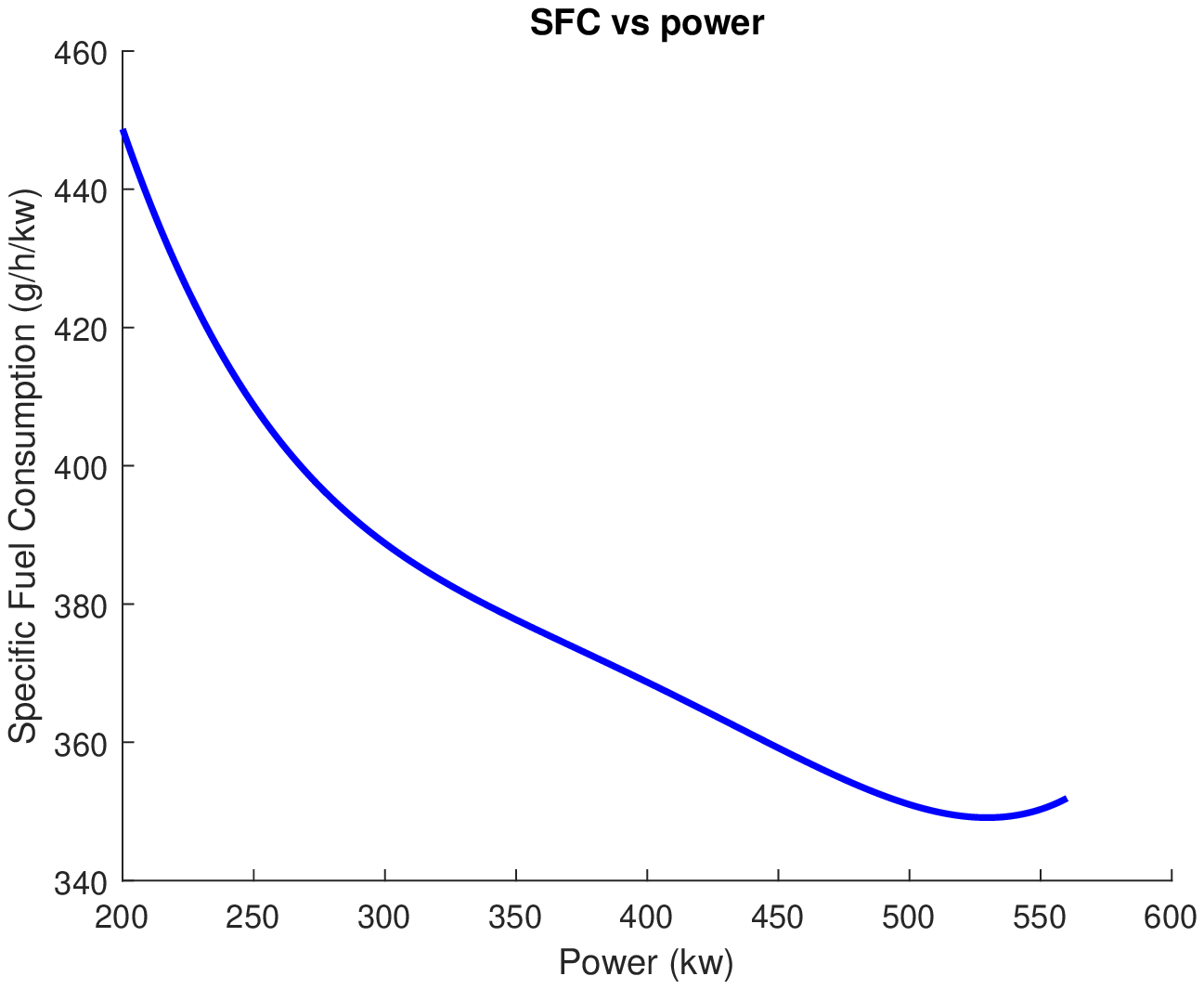}}
        \hfil
        \subfloat[]{\includegraphics[width=1.7in]{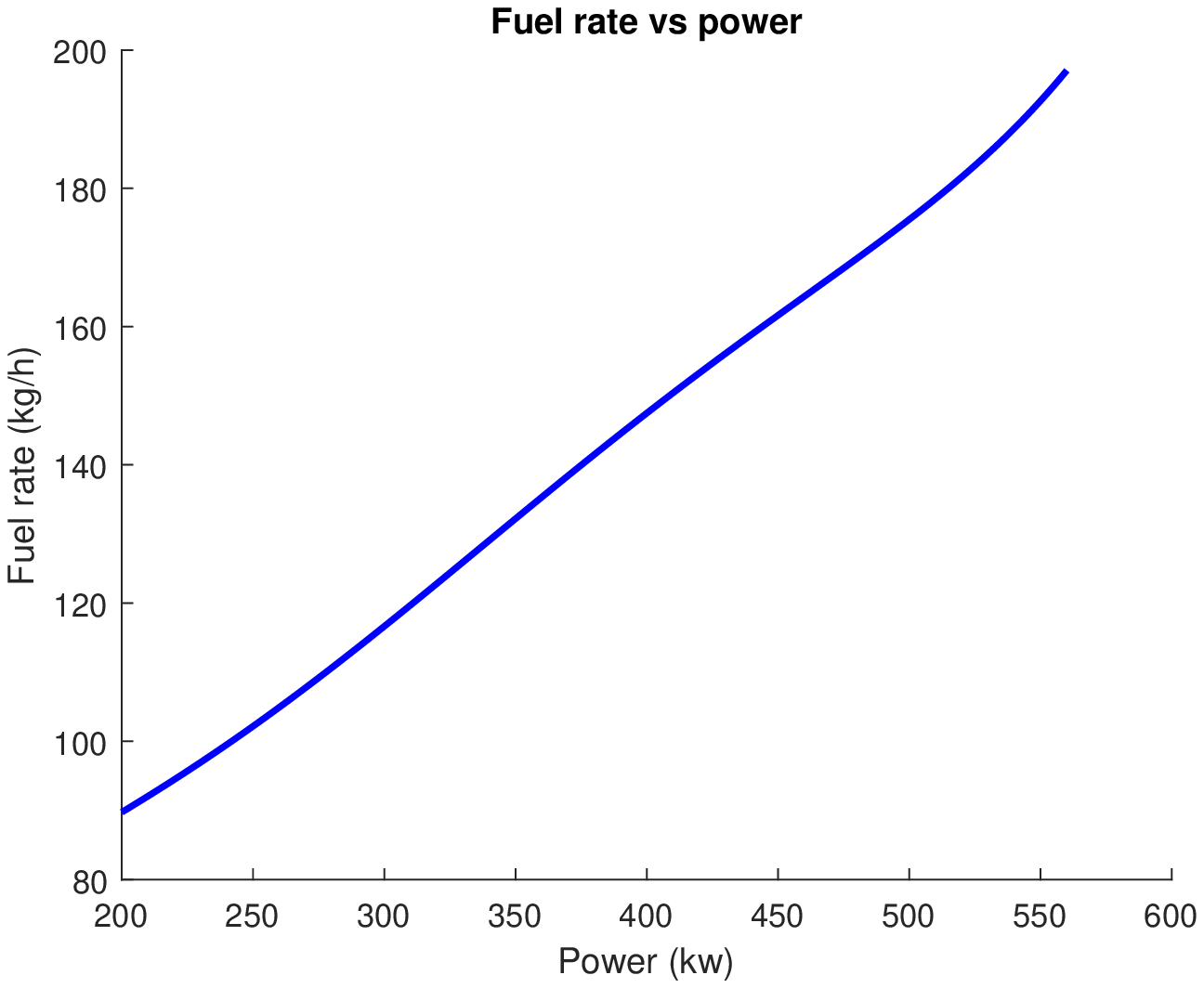}}
        \caption{Fuel consumption characteristics of the engine: (a)~SFC vs output power; (b)~fuel flow rate vs output power}
        \label{Fig_Preliminary_fuel}
    \end{figure}

    \item \textbf{Electric motor}: the Magni250 by MagniX is selected for our analysis; it weighs $72\mbox{kg}$; its maximum output power is $280\mbox{kW}$; it is assumed that the total efficiency of this motor is $95\%$;
    \item \textbf{Battery cells}: in~\cite{finger2020comparison}, the authors assumed a futuristic value for the energy density of the battery cell, significantly greater than the current state of battery technology. In this paper, an existing battery model from a Korean battery company Kokam is used for the analysis. The technical data of this battery is given in Table.~\ref{Table:NMC data}. 
    \begin{table}[h]
        \centering
        \caption{Technical data of an existing NMC cell}
        \label{Table:NMC data}
        \begin{tabular}{c|c|c|c}
        \hline
            Model & Capacity  & Weight & \makecell{Energy\\ Density} \\
            \hline
             SLPB140460330 & 200$\mbox{Ah}$ & 3.96 $\mbox{kg}$ & 189$\mbox{Wh/kg}$\\
             \hline
        \end{tabular}
    \end{table}
    For the battery pack configuration, it is assumed that all battery cells operate identically, and that the parallel and series configurations of the battery pack do not affect the output voltage and the state of charge of each cell. In addition, it is assumed that battery cell does not have an internal state dynamics, resulting in a constant output voltage of $3.6\mbox{V}$. For the power allocation algorithm in Section~\ref{Sec_PowerAllocation}, it is assumed that the battery pack has 2 parallel paths and 150 cells connected in series in each path. This configuration has an output voltage of $540\mbox{V}$ and a total capacity of $400\mbox{Ah}$. The total weight of the battery pack is $600\mbox{kg}$, which is greater than the value in~\cite{finger2020comparison}. Other battery pack configurations are also examined in Section~\ref{Sec:comparison}.
\end{enumerate}

\section{Flight dynamics and power relation in the propulsion system}\label{Sec_FlightDynamics}
Developing effective power allocation algorithms for HEA hinges upon interdependence between the aircraft flight dynamics and its propulsion system. In this section, we provide a brief overview of these two subsystems; the flight dynamics is introduced in Section~\ref{Subsec_Longitudinal dynamics} and the power relation in the aircraft propulsion system in Section~\ref{Subsec_propulsion system}.
\subsection{Longitudinal dynamics during the cruise phase}\label{Subsec_Longitudinal dynamics}
The complete longitudinal dynamics can be found in~\cite{wang2021power}; in this section, we focus on the optimal velocity profile during the cruise phase using the following assumptions: 
\begin{enumerate}
    \item The aircraft is flying from one point to another at a constant altitude;
    \item The angle of attack is small and as such $\sin{\alpha}\approx 0$ and $\cos{\alpha}\approx 1$.
\end{enumerate}
The simplified dynamics during the cruise phase is,
\begin{subequations}\label{Eqn_dyn_cruise}
    \begin{align}
        \dot{y}&=v\\
        \dot{v}&=\frac{1}{m}(T-D)\\
        L&=W,
    \end{align}
\end{subequations}
where $y$ is the horizontal position, $v$ is the velocity, $T$ is the thrust force, $D$ is the drag force, and $m$ is the total mass of the aircraft. Moreover, $W$ is the total weight of the aircraft and $L$ is the lift force. The drag and lift forces can be expressed as,
\begin{subequations}
    \begin{align}
        D&=\frac{1}{2}\rho SC_{\small D}v^2\\
        C_{\small D}&=C_{\small D,0}+KC_{\small L}^2\\
        L&=\frac{1}{2}\rho SC_{\small L}v^2,
    \end{align}
\end{subequations}
where $\rho$ is the air density; $S$ is the wing area; $C_{\small D}$ is the drag coefficient; $C_{\mbox{\small D,0}}$ is the zero lift drag coefficient; $K$ is the induced drag coefficient; all of these parameters are constant during the cruise phase. 

\subsection{HEA Propulsion system}\label{Subsec_propulsion system}
Then we discuss the power demand and power available in the aircraft propulsion system, which are required to subsequently formulate the power allocation problem.
\subsubsection{Power demand during the cruise phase}
The power demand $P_{\mbox{req}}$ during the cruise phase is 
\begin{equation}\label{Eqn_power_demand}
    P_{\mbox{\small req}}=Tv,
\end{equation}
\subsubsection{Power provided by the gas turbine}
We follow the assumption in~\cite{de2019preliminary} that the output power of the engine is given by
\begin{equation}\label{Eqn_P_e,out}
    P_{\mbox{\small e,out}}=\tau P_{\mbox{\small e,max}},
\end{equation}
where $\tau$ is the throttle of the engine, and $P_{\mbox{\small e,max}}$ is the maximum output power of the engine. 
\subsubsection{Power provide by the electric motor}
It is assumed that the transfer efficiency of the electric motor is a constant value; hence, the output power from the electric motor is given as,
\begin{equation}\label{Eqn_P_m,out}
    P_{\mbox{\small m,out}}=\eta_m P_{\mbox{\small bat}}=\eta_m UI ,
\end{equation}
where $\eta_m$ is the total efficiency of the electric motor; $P_{\mbox{\small bat}}$ is the output power of the battery pack; $U$ is the voltage of the battery pack; $I$ is the current passing through the battery pack.
\subsubsection{HEA Power relation}
The power relation in the propulsion system can be presented as
\begin{equation}\label{Eqn_power_relation}
   P_{\mbox{\small req}}=\eta_p(P_{\mbox{\small e,out}} + P_{\mbox{\small m,out}}) ,
\end{equation}
where $\eta_p$ is the transfer efficiency of the propeller, assumed to be a constant value for our analysis. Using Eq.(\ref{Eqn_power_demand})- Eq.(\ref{Eqn_power_relation}), we thus have the following useful relation,
\begin{equation}\label{Eqn_power_relation_expand}
    Tv=\eta_p(\tau \cdot P_{\mbox{\small e,max}}+\eta_m UI).
\end{equation}

\section{Power allocation algorithm for HEA}\label{Sec_PowerAllocation}
In this section, the HEA power allocation problem (originally represented as an optimal control problem) is formulated as a finite-dimensional optimization; then the second-order sufficient conditions for the solution of this problem are verified in Section~\ref{Subsec:theoretical analysis}. Numerical simulations supporting the theoretical analysis are then provided in Section~\ref{Subsec: numerical case}.
\subsection{Theoretical analysis}\label{Subsec:theoretical analysis}
In~\cite{wang2021power}, a power allocation algorithm for HEA, providing optimal power distribution between the engine and the electric motor for a given flight profile, has been proposed. In this context, it was assumed that the aircraft has a connected hybrid structure and the engine can charge the battery during flight. 
More specifically, the optimal control problem was formulated as,
\begin{subequations}\label{Eqn_power_allocation_ocp}
    \begin{align}
        J = &\min_{\tau}\quad -m(t_f)\\
        \mbox{s.t.}\quad & \dot{m}=-c\\
        & \dot{\theta}=-\frac{I}{Q}\\
        & m(t_0)=m_0; \quad \theta(t_0)=\theta_0\\
        & 0.1\leq \tau \leq 1\\
        & m_{\min}\leq m\leq m_{\max},\quad \theta_{\min}\leq \theta\leq \theta_{\max},
    \end{align}
\end{subequations}
where $\theta$ is the state of charge of the battery pack (SOC), and the current $I$ is a function of $m$ and $\tau$, deduced from Eq.~(\ref{Eqn_power_relation_expand}). In this paper, the flow rate $c$ is nearly a linear function of the output power of the engine. Consequently, instead of using $\mbox{SFC}$ and power to compute the flow rate, a representation of the flow rate as a function of the output power is derived and the dynamics of the aircraft mass is 
re-parameterized as a linear function of the control input as, 
\begin{equation}\label{Eqn_linear_dyn_m}
    \dot{m} = -(c_1P_{\mbox{e,out}}+c_2P_{\mbox{e,out}})=\tilde{c}_1\tau+\tilde{c}_2.
\end{equation}

In order to study the impact of the battery pack capacity on the fuel consumption, instead of the SOC of the battery pack, the remaining charge $q$ in the battery pack is taken as a state, whose dynamic is given by,
\begin{equation}\label{Eqn_new_state}
    \dot{q}=-I.
\end{equation}

The performance and solution of an optimal control problem are, in theory, invariant with respect to the magnitude of the states and the control input. In practice, however, numerical issues can arise in the solver when these magnitudes are highly different; in Problem~(\ref{Eqn_power_allocation_ocp}), the aircraft mass is approximately $6000\mbox{kg}$, and the range of the state of charge of the battery pack is 0 to 1. Hence, we use the standard remedy of applying the affine variable transformations to the aircraft mass and the remaining charge in the battery pack as,
\begin{subequations}\label{Eqn_scaled_states}
    \begin{align}
        \hat{m}&=a_1 m + b_1\\
        \hat{q}&=a_2 q + b_2;
    \end{align}
\end{subequations}
we now define the scaled time $\hat{t}\in[0,1]$; the dynamics of ``scaled states'' with respect to the ``scaled time'' is now written as,
\begin{subequations}\label{Eqn_scaled_states_dyn}
    \begin{align}
        \dot{\hat{m}}&=k_{11}\tau + k_{10}\\
        \dot{\hat{q}}&=k_{21}\tau + k_{20} + k_{22}\hat{m}^2 + k_{23}\hat{m},
    \end{align}
\end{subequations}
where both states are linear in the control input $\tau$. The scaled optimal control problem is now written as,
\begin{subequations}\label{Eqn_scaled_OCP}
    \begin{align}
        J = &\min_{\tau}\quad -\hat{m}(\hat{t}=1)\\
        \mbox{s.t.}\quad & \dot{\hat{m}}=k_{11}\tau + k_{10}\\
        & \dot{\hat{q}}=k_{21}\tau + k_{20} + k_{22}\hat{m}^2 + k_{23}\hat{m}\\
        & \hat{m}(0)=1; \quad \hat{q}(0)=\hat{q}_0\\
        & \hat{q}(1)=\hat{q}_f\\
        & 0\leq \tau \leq 1\\
        & 0\leq \hat{m}\leq 1,\quad 0\leq \hat{q}\leq 1.
    \end{align}
\end{subequations}

This formulation of the power allocation problem is in the form of bang-bang control problem with state constraints. When solving this problem using Tomlab for a reasonable flight profile, the only active state constraint is $\hat{q}=1$; hence, we subsequently only consider the constraint $\hat{q}\leq 1$. This constraint is denoted as,
\begin{equation}\label{Eqn_scaled_state_constr}
    S(\hat{q}(\hat{t})) = \hat{q}-1 \leq 0.
\end{equation}
 $S(\hat{q}(t))\leq 0$ is denoted as $S$ for simplicity in the sequel. An interval $[\hat{t}_1,\hat{t}_2]\in[0,1]$ is called a boundary arc if $S(\hat{q}(t))\equiv 0$ holds for all $t\in[\hat{t}_1,\hat{t}_2]$. Hence on the boundary arc the time derivative of $S$ with respect to the scaled time $\hat{t}$ should be zero for all $t\in[\hat{t}_1,\hat{t}_2]$, where,
\begin{equation}\label{Eqn_scaled_state_constr_dyn}
\frac{dS}{d\hat{t}}=\frac{\partial S}{\partial \boldsymbol{x}}\frac{\partial \boldsymbol{x}}{\partial \hat{t}} = [0~1]\begin{bmatrix}\dot{\hat{m}}\\ \dot{\hat{q}}\end{bmatrix}=\dot{\hat{q}},
\end{equation}
with $\boldsymbol{\mbox{x}}=[\hat{m} ~\hat{q}]^\top$.
Hence, the boundary control $\tau_b$ can be computed from $\frac{dS}{dt}=\dot{\hat{q}}=0$ as the following feedback expression
\begin{equation}\label{Eqn_bdn_arc_ctr}
    \tau_b=-\frac{1}{k_{21}}(k_{20} + k_{22}\hat{m}^2 + k_{23}\hat{m}).
\end{equation}

The augmented Hamiltonian of Problem (\ref{Eqn_scaled_OCP}) is given as,
\begin{align}\label{Eqn_scaled_Hamiltonian}
    H&=\lambda_1 \cdot \dot{\hat{m}} + \lambda_2 \cdot \dot{\hat{q}} + \eta \cdot S\nonumber\\
    &=\lambda_1 k_{10} + \lambda_2(k_{20} + k_{22}\hat{m}^2 + k_{23}\hat{m}) + \eta(\hat{q}-1)\nonumber\\
    &+ (\lambda_1 k_{11} + \lambda_2 k_{21})\tau \nonumber\\
    &=H_1 + \sigma (\hat{t})\tau.
\end{align}

Notice that the Hamiltonian is linear in the control $\tau$;
from the Pontryagin's minimum principle, the optimal control law on interior arcs with $S< 0$ is obtained as,
\begin{equation}\label{Eqn_scaled_ctr_strc}
    \tau^*=\Big \{ \begin{array}{cc}
         \tau_{\max}& \mbox{if} ~\sigma(t) <0 \\
         \tau_{\min} & \mbox{if} ~\sigma(t) >0.
    \end{array}
\end{equation}

For a boundary arc, the optimal control has been computed as $\tau_{\min}< \tau_b<\tau_{\max}$, since the optimal control should minimize the Hamiltonian, this minimum principle yields,
\begin{equation}\label{Eqn_bnd&singular}
    \sigma(\hat{t})=0,\quad \forall \hat{t}\in [\hat{t}_1,\hat{t}_2].
\end{equation}

This relation is interpreted as the property that a boundary control behaves formally like a singular control by Maurer in~\cite{maurer1977optimal}. When solving Problem (\ref{Eqn_scaled_OCP}) using a numerical solver, the optimal control follows a ``max-boundary-min'' structure. Combining the numerical results with the preceding theoretical analysis, the following optimal control structure is thus obtained,
\begin{equation}\label{Eqn_scaled_ctr_strc_withbnd}
    \tau^*=\Bigg \{ \begin{array}{cc}
         \tau_{\max}& 0\leq \hat{t} \leq \hat{t}_1 \\
         \tau_{b} & \hat{t}_1\leq \hat{t} \leq \hat{t}_2\\
         \tau_{\min} & \hat{t}_2\leq \hat{t} \leq 1.
    \end{array}
\end{equation}

In order to study the second order sufficient conditions and the sensitivity to initial conditions of states, this optimal control problem is then formulated as a (finite-dimensional) optimization problem. This approach is particular useful when the structure of the optimal control, i.e., the sequence of bang-bang and boundary arcs, has been derived. The \textit{arc-parametrization} method in~\cite{kaya2003computational} is applied to formulate the optimization problem with variables $\textbf{x}=[\xi_1,~\xi_2,~\xi_3,~z]^\top$(different $\boldsymbol{x}$ than in Eq.~(\ref{Eqn_scaled_state_constr_dyn})), as,
\begin{subequations}\label{Eqn_scaled_opt_variables}
    \begin{align*}
        \xi_1&=\hat{t}_1-0,          & \xi_2=\hat{t}_2-\hat{t}_1\\
        \xi_3&=\hat{t}_f-\hat{t}_2,  & z=\hat{m}(\hat{t}_2);
    \end{align*}
\end{subequations}
the last variable is introduced for the simplicity formulation of constraints. The power allocation problem can now be represented as an optimization problem,
\begin{subequations}\label{Eqn_scaled_opt}
    \begin{align}
        \min_{\tau}~ G\coloneqq &~-z - k_{10}\xi_3\\
        \mbox{s.t.}  ~\Phi_1\coloneqq &h_{13}\xi_1^3 + h_{12}\xi_1^2 + h_{11}\xi_1 + h_{10} - 1 = 0\\
         \Phi_2\coloneqq & g_{13}\xi_3^3+g_{12}z\xi_3^2 + g_{11}\xi_3^2+k_{20}\xi_3\nonumber\\
        &~ + k_{22}z^2\xi_3 + k_{23}z\xi_3 + 1 - \hat{q}_f = 0\\
         \Phi_3\coloneqq&~ z- \frac{-\tan(\frac{r_0\xi_2}{2})+f_1m_0+f_2\xi_1+f_3}{f_1(1+\tan(\frac{r_0\xi_2}{2})(f_1m_0+f_2\xi_1+f_3))}\nonumber\\
        & ~+ \frac{k_{31}}{2k_{32}} = 0\\
         \Phi_4\coloneqq &~\xi_1 + \xi_2 + \xi_3 =1;
    \end{align}
\end{subequations}
the detailed process of deducing this above optimization problem from the original optimal control problem is provided in Appendix \ref{App_trans2opt}. Note that the Lagrangian for Problem Eq.~(\ref{Eqn_scaled_opt}) is
\begin{equation}\label{Eqn_opt_Lagrangian}
    L \coloneqq G + \rho_1 \Phi_1 + \rho_2 \Phi_2 + \rho_3 \Phi_3 + \rho_4 \Phi_4,
\end{equation}
that is useful to derive second-order sufficient conditions in Appendix \ref{App_trans2opt}. 

\subsection{Numerical case study}\label{Subsec: numerical case}
\begin{figure}[h]
    \subfloat[]{\includegraphics[width=0.24\textwidth]{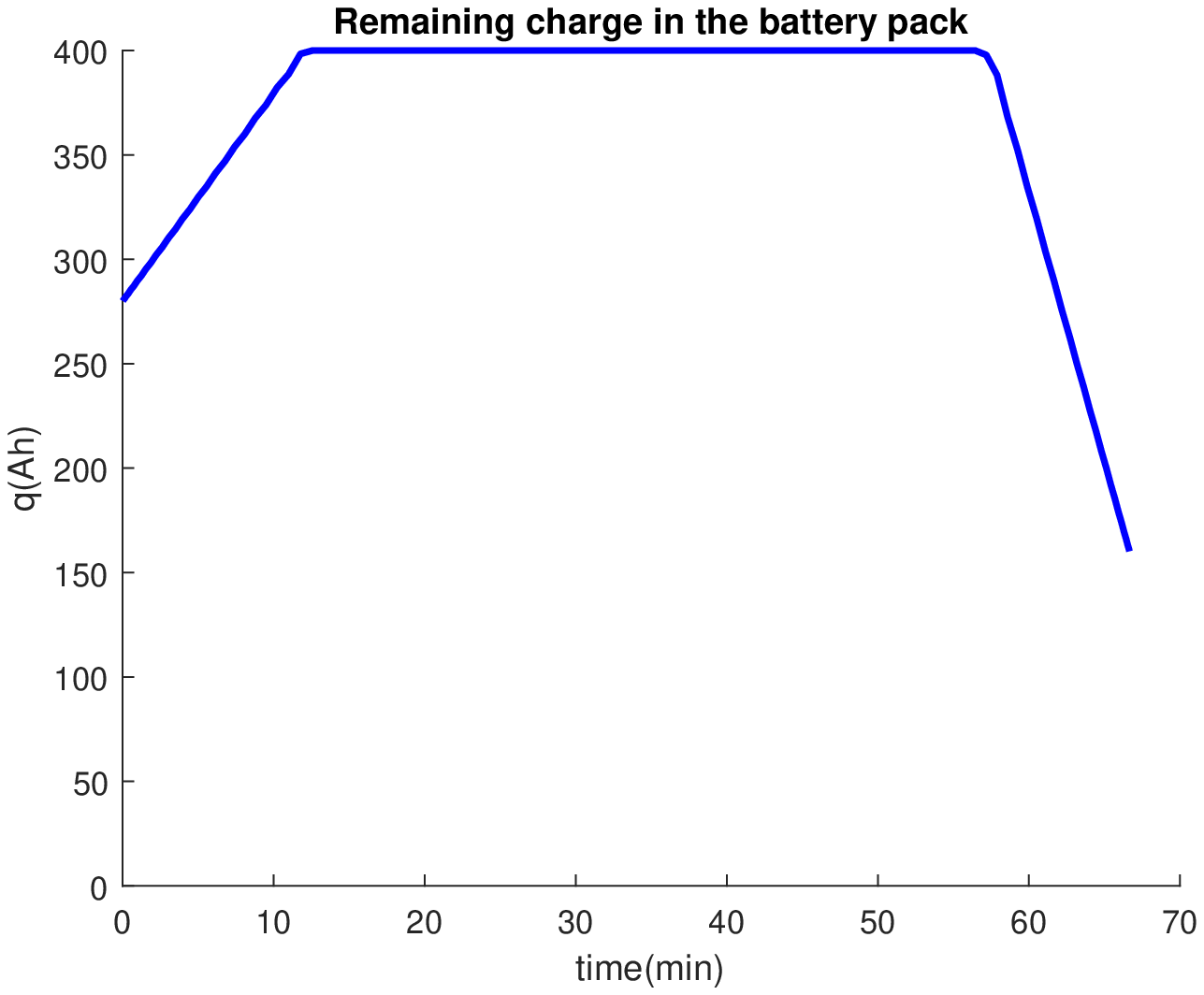}}
    \subfloat[]{\includegraphics[width=0.24\textwidth]{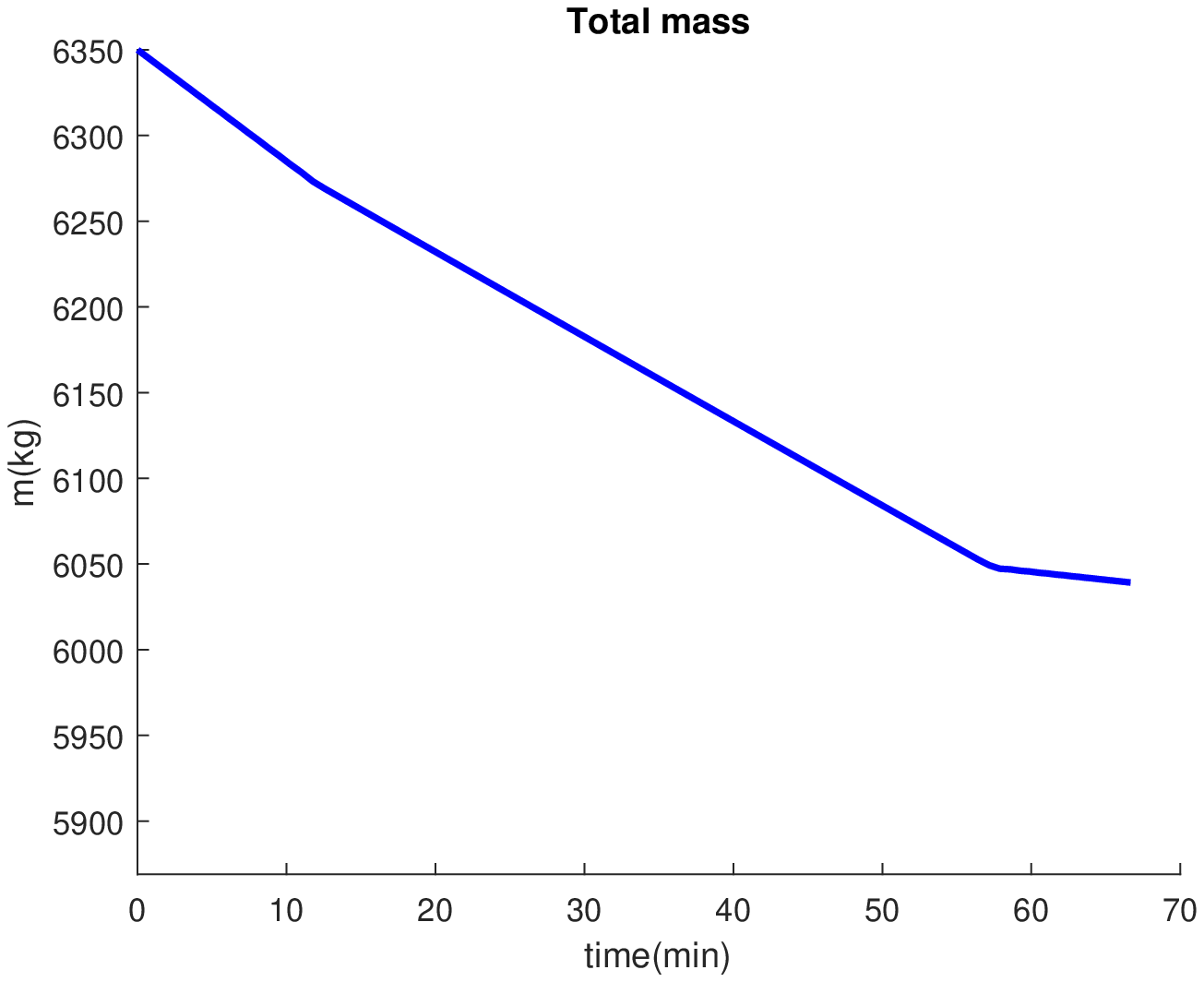}}
    \caption{(a) Remaining charge in the battery pack; (b)Total mass of the aircraft}
    \label{Fig_PA_SOC&mass}
\end{figure}

\begin{figure}[h]
    \subfloat[]{\includegraphics[width=0.24\textwidth]{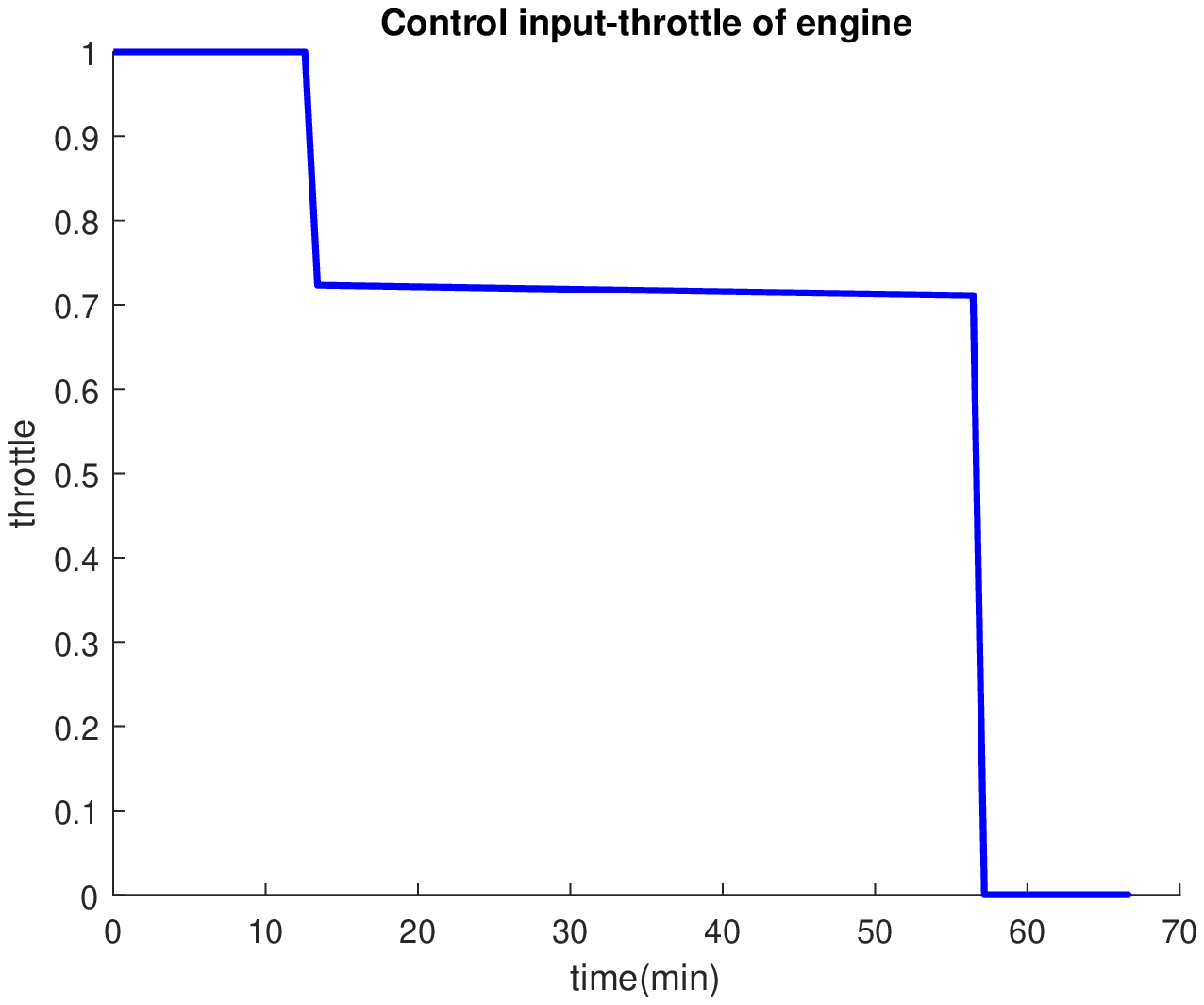}}
    \subfloat[]{\label{Fig_PA_sfc}\includegraphics[width=0.24\textwidth]{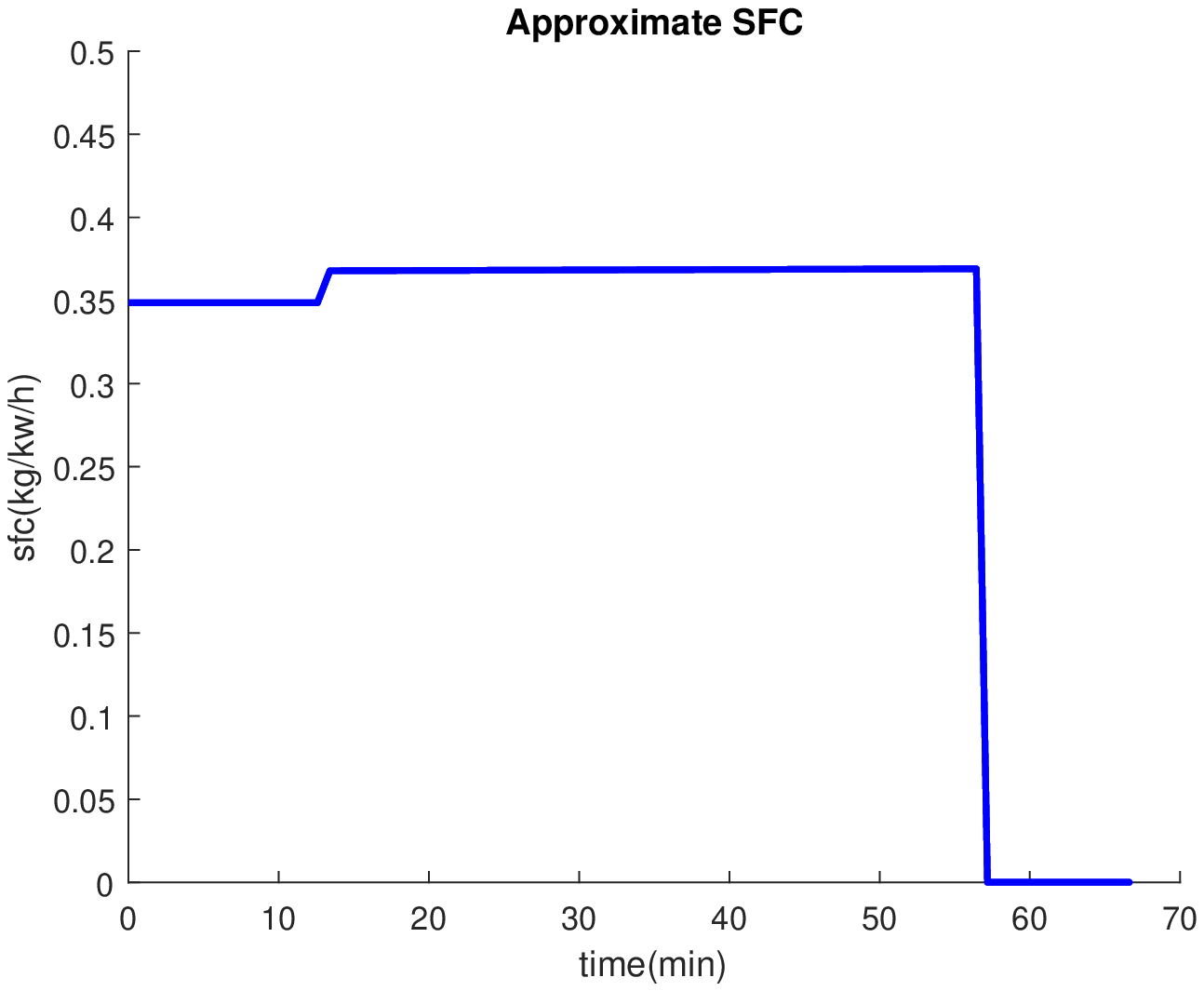}}
    \caption{(a) Control input; (b) Approximate Specific Fuel Consumption}
    \label{Fig_PA_control&sfc}
\end{figure}

We now consider the following flight conditions to examine the power allocation problem for the HEA: $y_d = 400 \mbox{km}$, $v=100\mbox{m/s}$, $\theta_0=0.7$, $m_0=6350\mbox{kg}$, $\theta_f = 0.4$; we use Tomlab to solve Problem (\ref{Eqn_scaled_OCP}) and re-scale the results back to original states. The results are given in Fig.~\ref{Fig_PA_SOC&mass} and Fig.\ref{Fig_PA_control&sfc}. The fuel consumed for this example is $\mbox{310.8kg}$. Fig.~\ref{Fig_PA_control&sfc} shows that the optimal control structure follows the optimal control structure in Eq.~(\ref{Eqn_scaled_ctr_strc_withbnd}).

Then the corresponding optimization problem is solved using ipopt; the optimal variables are determined as, $[\xi_1,~\xi_2,~\xi_3,~z]^\top = [0.1807, 0.6801, 0.1392, 0.3706]$, and the optimal objective value is $-\hat{m}(1) = -0.3525$. Hence the aircraft mass at the final time is $\mbox{m=6039.2kg}$, and the fuel consumed is $\mbox{310.8kg}$. In this numerical case, the dual variables are also provided by ipopt, namely, $[\rho_1, \rho_2,\rho_3,\rho_4]=[-0.1253,-0.1272,0.9971,-0.6783]$; the first order necessary condition for optimality can now be checked, i.e., $L_{\boldsymbol{x}}=0$; moreover, the Hessian $L_{\boldsymbol{xx}}$ of the Lagrangian for $\boldsymbol{x}=[0.1807, 0.6801, 0.1392, 0.3706]$ is 
\begin{equation}
    L_{\boldsymbol{xx}}=\begin{bmatrix}
    -0.0186 & -0.0188 & 0 & 0\\
    -0.0188 & -0.0144 & 0 & 0\\
    0 & 0 & -0.0028 & 0.0211\\
    0 & 0 &0.0211 & 0.0002,
    \end{bmatrix}
\end{equation} 
and 
\begin{equation}
    \Phi_{\boldsymbol{x}}=\begin{bmatrix}
    1.6741 & 0 & 0 & 0\\
    0 & 0 & -4.3085 & -0.0231\\
    0.8907 & 0.6803 & 0 & 1\\
    1 & 1 & 1 & 0
    \end{bmatrix},
\end{equation}
with $\mbox{\bf rank}\,\Phi_{\boldsymbol{x}}=4$ and $L_{\boldsymbol{xx}}$ is positive definite. Hence, the second-order sufficient conditions are satisfied, and the optimal solution is unique and, in fact, globally optimal. With this transition to the finite-dimensional optimization problem, one can conduct sensitivity analysis on the optimal solution, a technique that proves instrumental to address the question posed in the title of this manuscript.
\section{To Charge in-flight or not: comparison of two parallel hybrid architectures}\label{Sec:comparison}
The mechanical connection between the engine and the electric motor in the HEA propulsion system allows the engine to charge the battery pack in-flight. For the climb phase, the power demand is usually higher than the maximum output power of the engine--hence the propulsion system can only run in the combined mode. The power-split mode is only applicable during the cruise phase and the descent phase. In this section, it is demonstrated that the fuel saving by charging the battery during the cruise phase to steer the engine to operate more efficiently is negligible. The second potential advantage of charging the battery in-flight is that the electric motor can work in OEI mode by charging the battery to its fully charged state after the climb phase. It has been shown that under certain flight conditions, this advantage can be achieved instead by carrying more battery packs before take-off.

\subsection{To Charge or not in-flight during the cruise phase}
First, fuel consumption between charging and not charging the battery pack in the cruise phase is compared. By solving Problem (\ref{Eqn_scaled_OCP}) with an extra constraint $\dot{\hat{q}}\leq 0$, one can determine the optimal fuel consumption when charging the battery pack during the flight is not allowed (the current passing through the battery pack can only be positive). In
this direction, it is assumed that all other conditions and constraints remain the same. 
\begin{figure}[h]
    \subfloat[]{\includegraphics[width=0.24\textwidth]{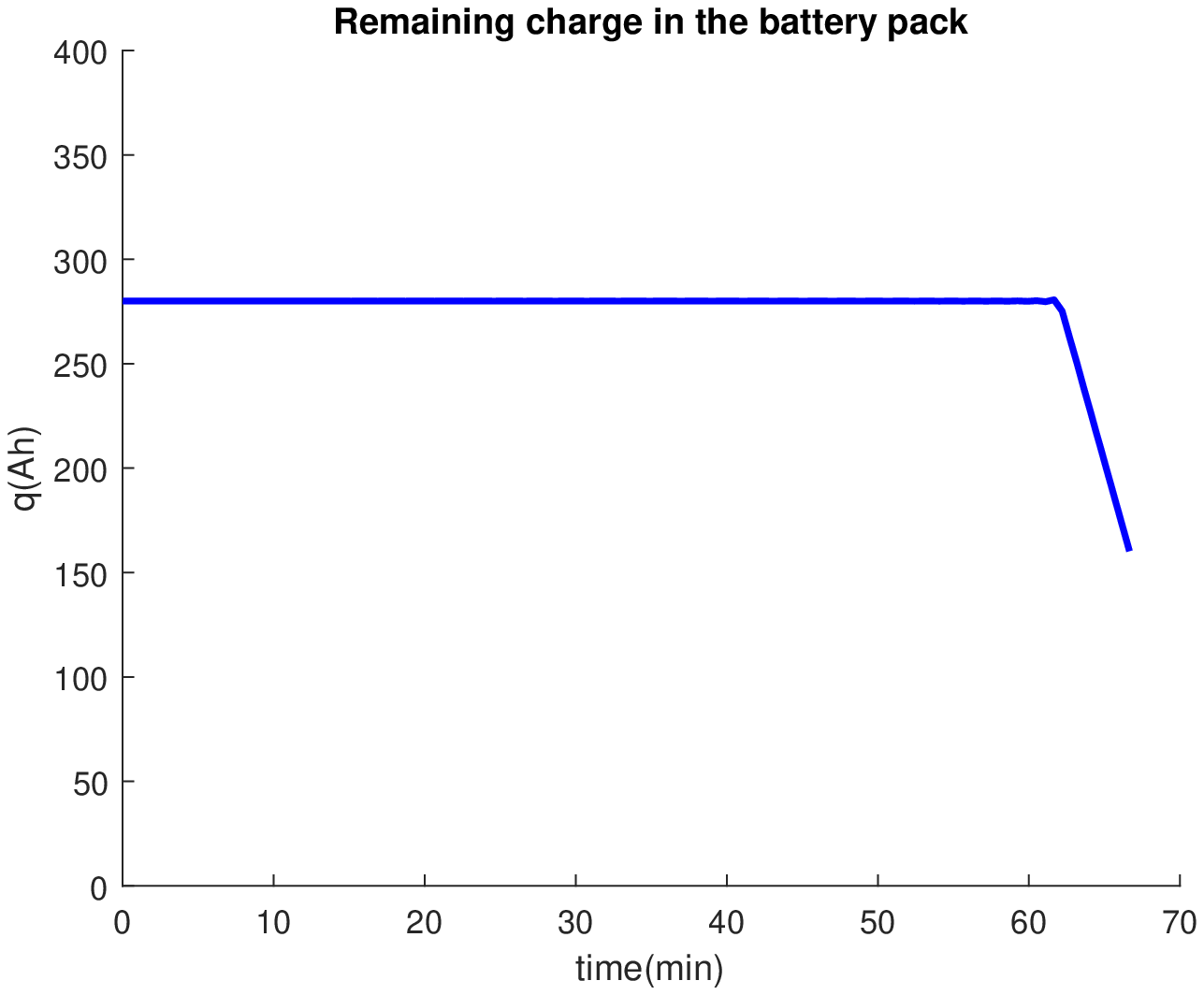}}
    \subfloat[]{\label{Fig_PA_sfc_notcharge}\includegraphics[width=0.24\textwidth]{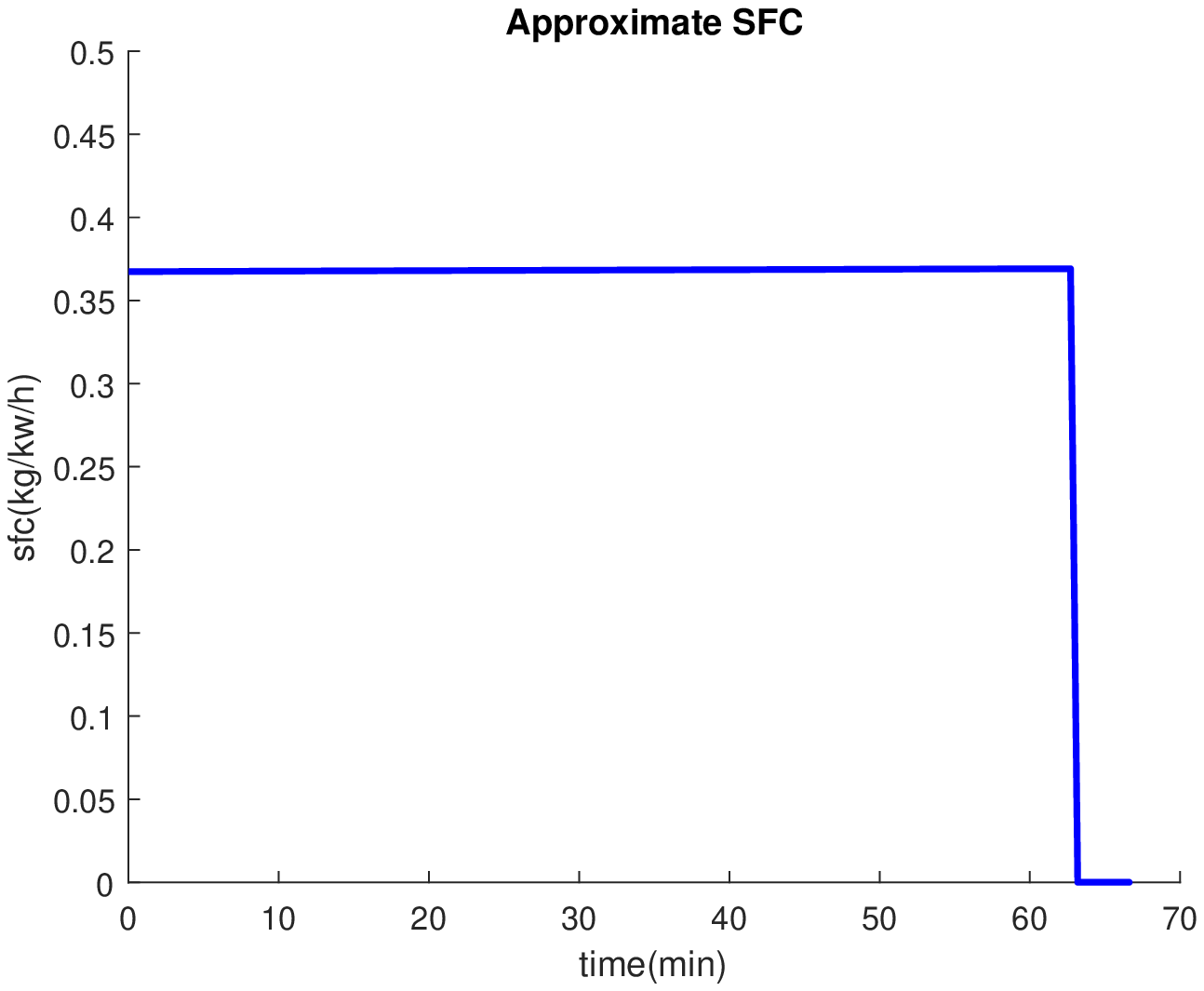}}
    \caption{(a) Remaining charge in the battery pack; (b) Approximate Specific Fuel Consumption}
    \label{Fig_PA_soc&sfc_notcharge}
\end{figure}
Observing Fig.~\ref{Fig_PA_soc&sfc_notcharge}, when not allowing charging of the battery during the flight, the optimal power distribution strategy yields an ``engine-alone'' then ``motor-alone'' structure. The total fuel consumed is $311.08\mbox{kg}$, which is only $0.33\mbox{kg}$ more than the case when the battery was allowed to charge in-flight. Comparing the first segment of Fig.~\ref{Fig_PA_sfc} and Fig.~\ref{Fig_PA_sfc_notcharge}, we note that when charging the battery on $[0,\hat{t}_1]$ in Fig.~\ref{Fig_PA_sfc}, the SFC is about $0.35\mbox{kg/kw/h}$, which is lower than the case when battery was not charging ($0.367-0.369\mbox{kg/kw/h}$) in Fig.\ref{Fig_PA_sfc_notcharge}. However after charging the battery to its full state, the aircraft will run in the engine alone mode again, where the SFC is about $0.367-0.369\mbox{kg/kw/h}$ in Fig.~\ref{Fig_PA_sfc}, same as the case when charging was not not permissible. 

Charging the battery during flight can make the engine work more efficiently, but the duration of this efficient operation depends on the initial and terminal constraints of the battery pack. After charging the battery to its full capacity, the engine will go back to its inefficient operational region. In fact, overall, the fuel saving by charging the battery to steer the engine operation more efficiently is rather negligible.

\subsection{Carrying more battery in-flight}
The second benefit of charging the battery is that the electric motor can work in the OEI mode. In order to verify this, we compare the total fuel consumption for the combination of the climb and the cruise phase when the take-off weight is different due to different battery pack configurations. In the previous section, it is assumed that the battery pack has 2 parallel paths with a capacity of $400\mbox{Ah}$. Each parallel path in the battery pack has a capacity of $400\mbox{Ah}$ and weighs $600\mbox{kg}$. Different battery pack configurations are achieved by changing the number of parallel paths. Adding/removing parallel paths only increases/decreases the total capacity, not affects the output voltage of the battery pack. First, the fuel consumption during the climb phase with a fixed flight path angle is computed, and the terminal condition for the climb phase is when the altitude reaches $3000\mbox{m}$. 
\begin{table}[h]
    \centering
    \caption{Fuel consumption of the climb phase with different battery pack configuration}
    \begin{tabular}{c|c|c|c}
    \hline
         Number of parallel paths & 1 & 2 & 3  \\
         \hline
         \makecell{Total capacity of\\ the battery pack ($\mbox{Ah}$)}
          & 200 & 400 & 600 \\
         \hline
         MTOW ($\mbox{kg}$) & 5785 & 6385& 6985 \\
         \hline
         SOC after climb & 0.38 & 0.69 & 0.79 \\
         \hline
         Remaining charge ($\mbox{Ah}$) & 76 & 276 & 474 \\
         \hline
         Fuel consumed during climbing ($\mbox{kg}$) & 32.84 & 37.55 & 42.02 \\
         \hline
    \end{tabular}
    \label{Table:fuel climb}
\end{table}

From Table.~\ref{Table:fuel climb}, it is observed that the consumed fuel difference is only $10\mbox{kg}$ when the take-off weight difference is $1200\mbox{kg}$ (between 1 path and 3 paths). The optimal power allocation solutions show that during the climb phase, the electric motor operates at its maximum output power, and the engine runs based on the power demand.

For the cruise phase, we compare the following cases:
\begin{enumerate}
    \item case 1: $Q=200\mbox{Ah}$, $m_0=5750\mbox{kg}$, $\theta_0=0.4$, $\theta_f=1$;
    \item case 2: $Q=400\mbox{Ah}$, $m_0=6350\mbox{kg}$, $\theta_0=0.7$, $\theta_f=1$;
    \item case 3: $Q=400\mbox{Ah}$, $m_0=6350\mbox{kg}$, $\theta_0=0.7$, $\theta_f=0.5$, and the engine cannot charge the battery;
    \item case 4: $Q=600\mbox{Ah}$, $m_0=6950\mbox{kg}$, $\theta_0=0.8$, $\theta_f=0.67$, and the engine cannot charge the battery.
\end{enumerate}

\begin{figure}
    \centering
    \includegraphics[width=0.45\textwidth]{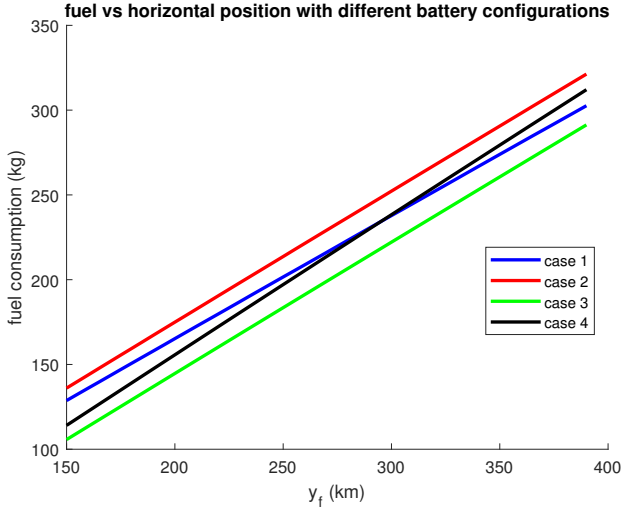}
    \caption{Fuel consumption with different battery configurations; case 1 and case 2 are connected configurations; case 3 and case 4 are independent configurations}
    \label{Fig_comp_diffQ}
\end{figure}

It is assumed that case 1 has the connected architecture; hence the engine can charge the battery in-flight. The terminal condition for the battery is to charge it to $200\mbox{Ah}$. The same architecture and terminal conditions are assumed for case 2. Cases 3 and 4 have independent architectures, the terminal condition is that the battery pack has $200\mbox{Ah}$ left. Then the cruise range from $150\mbox{km}$ to $400\mbox{km}$ is swept to compute the fuel consumption under the four cases listed above. The cruise speed is fixed as $80\mbox{m/s}$ for all cases. The simulation result is given in Fig.~\ref{Fig_comp_diffQ}.  

It is observed that the fuel consumption is nearly a linear function of the flight distance for each case (with different slope for each case). The sensitivity derivative shows that the derivative of the objective function with respect to the terminal time $\hat{t}_f$ is a constant value, but this value changes slightly as the terminal time changes. For example, in case 2, $\frac{dG}{dt_f}=0.6498$ when $t_f=0.375$; $\frac{dG}{dt_f}=0.6363$ when $t_f=1$; this value decreases as the final time $t_f$ increases. In the simulation result, the slope of the red line in Fig.~\ref{Fig_comp_diffQ} is 0.6433. Reference~\cite{kim2003sensitivity} provides explicit formula to compute the sensitivity derivatives,
\begin{equation}
    \begin{bmatrix}
    \frac{d \boldsymbol{x}^p}{dp}\\
    \frac{d \rho^p}{dp}
    \end{bmatrix}=-\begin{bmatrix}
    L_{\boldsymbol{xx}} & \Phi_{\boldsymbol{x}}^*\\
    \Phi_{\boldsymbol{x}} & \boldsymbol{0}
    \end{bmatrix}\begin{bmatrix}
    L_{\boldsymbol{x}p}\\
    \Phi_p
    \end{bmatrix},
\end{equation}
where $\boldsymbol{x}^p$ and $\boldsymbol{\rho}^p$ are the optimal and dual variables corresponding to parameter $p$. Then the derivative of the objective function $G$ with respect to the parameter $t_f$ is computed as
\begin{equation}
    \frac{d \mbox{G}}{d t_f}=\frac{\partial F}{\partial \boldsymbol{x}}\frac{\partial \boldsymbol{x}}{\partial t_f}.
\end{equation}

The following observations are obtained when comparing cases with the same terminal battery charge condition (after charging for the connected case; beginning of the cruise phase for the independent case):
\begin{enumerate}
    \item Comparing case 2 and case 3, it is observed that these two cases have the same initial aircraft mass and battery pack capacity; the only difference is the terminal condition for the SOC of the battery pack. Case 2 requires the battery to remain in a fully-charged state ($400\mbox{Ah}$) after a short period of the cruise flight; case 3 requires the remaining charge in the battery pack is more than $200\mbox{Ah}$. Simulation results show case 2 needs about $30\mbox{kg}$ more fuel consumed than case 3 for all flight ranges.
    \item Comparing case 1 and case 3 (the blue line and the green line), the remaining charge in the battery pack is $200\mbox{Ah}$ for both cases. Even though case 1 has lighter initial aircraft mass, charging the battery to its fully-charged state requires more fuel; carrying more battery pack can save fuel under this remaining charge condition.
    \item Comparing case 2 and case 4, same as the comparison between case 1 and case 3; the remaining charge in the battery pack is $400\mbox{Ah}$ for both cases. Carrying more battery can still save fuel.
\end{enumerate}

Hence, with the same remaining battery charge requirement, carrying more battery before take-off is more favorable than charging the battery in-flight in terms of fuel saving and $\mbox{CO}_2$ emissions. 

\section{Conclusions}\label{Sec:conclusion}
Two parallel hybrid electric architectures are compared in terms of fuel consumption. The connected architecture has the mechanical connection between the engine and the electric motor, allowing the engine to charge the battery in-flight in order to improve the engine's fuel efficiency. The independent architecture does not have a mechanical connection, which allows an easy transition from hybrid aircraft to all-electric aircraft. An optimal power allocation algorithm is applied to compare the fuel usage of these two architectures. Both theoretical and simulation results show that the fuel saving by charging the battery is negligible for the 19-seat conceptual aircraft. The other potential advantage of charging the battery during flight is also examined. The simulation results show that with the same requirement on the remaining battery charge, carrying more battery before taking-off is more fuel efficient than charging the battery during the flight. As such, the connected architecture cannot bring significant fuel savings compared with the independent configuration. In addition, the mechanical connection increases the system complexity and generally requires sophisticated control strategies. This architecture also makes it difficult to transition to an all-electric aircraft configuration. Lastly, charging the battery in the harsh flight environment requires accurate thermal management system as well as causing battery degradation. Our assessment--in the context of the question posed in the title of the manuscript-- is therefore implicit from these observations.

\section*{Appendix}
In this section, we provide details on some of the analytic steps used for our analysis.
\subsection{Transcription to an optimization problem}\label{App_trans2opt}
The explicit functions of the optimal control on each arc
are as follows:
\begin{enumerate}
    \item On $\hat{t}\in [0,\hat{t}_1]$, the optimal control is $\tau^*=\tau_{\max}=1$. Then the dynamics of states can be written as,
    \begin{subequations}\label{Eqn_dyn_max_arc}
        \begin{align}
            \dot{\hat{m}}&=k_{11}+k_{10}\\
            \dot{\hat{q}}&=k_{21}+ k_{20} + k_{22}\hat{m}^2 + k_{23}\hat{m}.
        \end{align}
    \end{subequations}
    We therefore obtain the explicit state function of time as,
    \begin{subequations}\label{Eqn_sol_max_arc}
        \begin{align}
            \hat{m}(\hat{t})&=m_{0} + (k_{11}+k_{10})\hat{t}\\
            \hat{q}(\hat{t})&=h_{13}\hat{t}^3 + h_{12}\hat{t}^2 +h_{11}\hat{t} + h_{10} , \label{Eqn_sol_max_arc_q}
        \end{align}
    \end{subequations}
    where,
    \begin{subequations}
        \begin{align*}
            h_{13}&=\frac{1}{3}k_{22}(k_{11}+k_{10})^2\\
            h_{12}&=(k_{22}m_0+\frac{1}{2}k_{23})*(k_{11}+k_{10})\\
            h_{11}&=k_{21}+k_{20}+k_{22}m_0^2+k_{23}m_0,
        \end{align*}
    \end{subequations}
    and $h_{10}=\hat{q}_0$ is computed from $\hat{q}(\hat{t}=0)=\hat{q}_0$.
    \item On $\hat{t}\in [\hat{t}_1,\hat{t}_2]$, $\hat{q}(\hat{t})= 0$ holds for all $\hat{t}$, and the optimal control is $\tau_b=-\frac{1}{k_{21}}(k_{20} + k_{22}\hat{m}^2 + k_{23}\hat{m})$. Hence, the dynamics of the aircraft mass can be written as,
    \begin{equation}
        \dot{\hat{m}}=-\frac{k_{11}}{k_{21}}(k_{20} + k_{22}\hat{m}^2 + k_{23}\hat{m})+k_{10};
    \end{equation}
    the initial condition for this interval is $\hat{m}(\hat{t}_1)=m_{0} + (k_{11}+k_{10})\hat{t}_1$. Thereby, by separation of variables, we obtain the explicit function of the aircraft mass on the boundary arc as,
    \begin{equation}
        \hat{m}(\hat{t})=\frac{1}{2k_{32}}[r_0 \tan(-\frac{r_0(\hat{t}-\hat{t}_1)}{2}+C_0)-k_{31}] ,
    \end{equation}
    where $k_{30}=\frac{k_{11}k_{20}}{k_{21}}-k_{10}$ $k_{31}=\frac{k_{11}k_{23}}{k_{21}}$, $k_{32}=\frac{k_{11}k_{22}}{k_{21}}$, and \begin{align*}
        C_0&=\frac{-2\arctan(\frac{2k_{32}m(\hat{t}_1)+k_{31}}{r_0})}{r_0}\\
        r_0&=\sqrt{4k_{32}k_{30}-k_{31}^2} .
    \end{align*}
    We denote $\hat{m}(t_2)$ as $z$, which is the third constraint in Problem Eq.~(\ref{Eqn_scaled_opt}).
    
    \item On $\hat{t}\in [\hat{t}_2, 1]$, the optimal control is $\tau^*=\tau_{min}=0$, 
    then the dynamics of states can be written
    \begin{subequations}
        \begin{align}
            \dot{\hat{m}}&=k_{10}\\
            \dot{\hat{q}}&=k_{20} + k_{22}\hat{m}^2 + k_{23}\hat{m}
        \end{align}
    \end{subequations}
    where the initial states on this arc are $[z, 1]$; the explicit functions of the states are, 
    \begin{subequations}
        \begin{align}
            \hat{m}(\hat{t}) &= \hat{m}(\hat{t}_2) + k_{10}(\hat{t}-\hat{t}_2)\\
            \hat{q}(\hat{t}) &= g_{13}(\hat{t}-\hat{t}_2)^3+g_{12}z(\hat{t}-\hat{t}_2)^2 + g_{11}(\hat{t}-\hat{t}_2)^2\nonumber\\
            & +k_{20}(\hat{t}-\hat{t}_2)+ k_{22}z^2(\hat{t}-\hat{t}_2) + k_{23}z(\hat{t}-\hat{t}_2) + 1 .\label{Eqn_dyn_min_arc_q}
        \end{align}
    \end{subequations}
\end{enumerate}
The first constraint for the optimization problem is when $q$ hits the boundary at the junction time $\hat{t}=\hat{t}_1$, given as
\begin{equation}
    \Phi_1 \coloneqq \hat{q}(\hat{t}_1)=1.
\end{equation}
When we substitute this constraint into Eq.~(\ref{Eqn_sol_max_arc_q}), we obtain,
\begin{equation}
    h_{13}\xi_1^3 + h_{12}\xi_1^2 + h_{11}\xi_1 + h_{10} =1
\end{equation}
The second constraint for the optimization problem is the terminal condition for $\hat{q}$, $\hat{q}(1)=\hat{q}_f$, after substituting into Eq.~(\ref{Eqn_dyn_min_arc_q}), we obtain 
\begin{equation}
    g_{13}\xi_3^3+g_{12}z\xi_3^2 + g_{11}\xi_3^2+k_{20}\xi_3
        + k_{22}z^2\xi_3 + k_{23}z\xi_3 + 1 - \hat{q}_f = 0.
\end{equation}

\subsection{Verification of SSC and sensitivity analysis}\label{App_ssc}
We refer Theorem 3.1 in \cite{maurer2005optimization} pertaining to the second-order sufficient condition for the optimization problem Eq.~(\ref{Eqn_scaled_opt}). Let $\boldsymbol{x}$ be a feasible solution for Problem Eq.~(\ref{Eqn_scaled_opt}). Suppose the following conditions hold: 
\begin{enumerate}
    \item $L_{\boldsymbol{x}}=0$ (first order necessary condition).
    \item $\mbox{rank}[\Phi_{\boldsymbol{x}} ]=r$, where $\mbox{r}$ is the number of constraints
    \item $v^\top L_{\boldsymbol{xx}}v >0, \forall v\in \mathbb{R}^{n}, v\neq 0, \Phi_{\boldsymbol{x}}v=0$, where $n$ is the number of variables in Problem Eq.~(\ref{Eqn_scaled_opt}).
\end{enumerate}

Let $\Delta = 1 + \tan(\frac{r_0}{2}\xi_2)(f_1x_{10}+f_2\xi_1+f_3)$ and $f_4 = f_1x_{10}+f_2\xi_1+f_3$.

The $L_{\boldsymbol{x}}$ is now computed as,
\begin{subequations}
    \begin{align}
        L_{\xi_1} &= \rho_1(3h_{13}\xi_1^2 + 2h_{12}\xi_1 + h_{11}) + \rho_4 \nonumber\\
        &+ \frac{-\rho_3f_2}{f_1} \frac{(1+\tan^2(\frac{r_0}{2}\xi_2))}{\Delta^2}\\
        L_{\xi_2} &=\frac{\rho_3 r_0}{2f_1}\frac{(1+f_4^2)\sec^2(\frac{r_0}{2}\xi_2)}{\Delta^2} + \rho_4\\
        L_{\xi_3}&=-k_{10} + \rho_4 + \rho_2(3g_{13}\xi_3^2 + 2g_{12}z\xi_3\nonumber \\
        &+ 2g_{11}\xi_3 + k_{20}+k_{22}z^2+k_{23}z) \\
        L_{z} &= -1 + \rho_3 + \rho_2(g_{12}\xi_3^2+2k_{22}\xi_3 z + k_{23}\xi_3) .
    \end{align}
\end{subequations}

The $L_{\boldsymbol{xx}}$ is represented as,
\begin{align}
        L_{\boldsymbol{xx}}&=\begin{bmatrix}
    \frac{\partial L_{\xi_1}}{\partial \xi_1} & \frac{\partial L_{\xi_1}}{\partial \xi_2} & \frac{\partial L_{\xi_1}}{\partial \xi_3} & \frac{\partial L_{\xi_1}}{\partial z}\\
    \frac{\partial L_{\xi_2}}{\partial \xi_1} & \frac{\partial L_{\xi_2}}{\partial \xi_2} & \frac{\partial L_{\xi_2}}{\partial \xi_3} & \frac{\partial L_{\xi_2}}{\partial z}\\
    \frac{\partial L_{\xi_3}}{\partial \xi_1} & \frac{\partial L_{\xi_3}}{\partial \xi_2} & \frac{\partial L_{\xi_3}}{\partial \xi_3} & \frac{\partial L_{\xi_3}}{\partial z}\\
    \frac{\partial L_z}{\partial \xi_1} & \frac{\partial L_z}{\partial \xi_2} & \frac{\partial L_z}{\partial \xi_3} & \frac{\partial L_z}{\partial z}
    \end{bmatrix}\nonumber \\
    &= \begin{bmatrix}
    L_{11} & L_{12} & 0 & 0\\
    L_{21} & L_{22} & 0 & 0\\
    0 & 0 & L_{33} & L_{34}\\
    0 & 0 &L_{43} & L_{44}
    \end{bmatrix},
\end{align}
where
\begin{subequations}
    \begin{align}
        L_{11}&=\rho_1(6h_{13}\xi_1 + 2h_{12}) \nonumber\\
        &+ \frac{2\rho_3f_2^2}{f_1}\frac{(\tan(\frac{r_0}{2}\xi_2)+\tan^3(\frac{r_0}{2}\xi_2))}{\Delta^3}\\
        L_{12}&=\frac{-\rho_3r_0f_2}{f_1}\frac{\sec^2(\frac{r_0}{2}\xi_2)(\tan(\frac{r_0}{2}\xi_2)-f_4)}{\Delta^3}\\
        L_{21} & =L_{12}\\
        L_{22} &= \frac{\rho_3r_0^2(1+f_4^2)\sec^2(\frac{r_0}{2}\xi_2)}{2f_1}\nonumber \\
        &\frac{(\tan(\frac{r_0}{2}\xi_2)(f_4\tan(\frac{r_0}{2}\xi_2)+1)-f_4\sec^2(\frac{r_0}{2}\xi_2))}{\Delta^3} .
    \end{align}
\end{subequations}

The function $\Phi_{\boldsymbol{x}}$ is now computed as,
\begin{subequations}
    \begin{align}
        \frac{\partial \Phi_1}{\partial \xi_1} &= 3h_{13}\xi_1^2 + 2h_{12}\xi_1 + h_{11},\quad \frac{\partial \Phi_1}{\partial \xi_2} = 0 \nonumber\\
        \frac{\partial \Phi_1}{\partial \xi_3} &= 0, \quad \frac{\partial \Phi_1}{\partial z}=0;\\
        \frac{\partial \Phi_2}{\partial \xi_1} & =0, \quad \frac{\partial \Phi_2}{\partial \xi_2} =0\nonumber\\ 
        \frac{\partial \Phi_2}{\partial \xi_3} &= 3g_{13}\xi_3^2+2g_{12}z\xi_3+2g_{11}\xi_3+k_{20}+k_{22}z^2 + k_{23}z\nonumber\\
        \frac{\partial \Phi_2}{\partial z}&=g_{12}\xi_3^2+2k_{22}\xi_3z + k_{23}\xi_3;\\
        \frac{\partial \Phi_3}{\partial \xi_1} & =\frac{-f_2(1+\tan^2(\frac{r_0}{2}\xi_2))}{f_1\Delta^2},
        ~\frac{\partial \Phi_3}{\partial \xi_2}  = \frac{r_0\sec^2(\frac{r_0}{2}\xi_2)(1+f_4^2)}{2f_1\Delta^2}\nonumber\\
        \frac{\partial \Phi_3}{\partial \xi_3}&=0,\quad \frac{\partial \Phi_3}{\partial z}=1;\\
        \frac{\partial \Phi_4}{\partial \xi_1} &=1, \quad \frac{\partial \Phi_4}{\partial \xi_2} =1, \quad \frac{\partial \Phi_4}{\partial \xi_3} =1,\quad \frac{\partial \Phi_4}{\partial z} =0,
    \end{align}
\end{subequations}

\section*{ACKNOWLEDGMENT} The authors would like to thank Kuang-Ying Ting, Yue Yu, Professor Mohammad Reza Soltani, and Professor Krishna A. Shah for helpful discussions.

\vfill


\begin{thebibliography}{}
\bibitem{ye2021review} Y. Xie, A. Savvarisal, A. Tsourdos, D. Zhang, J. GU, ``Review of hybrid electric powered aircraft, its conceptual design and energy management methodologies,'' Chinese Journal of Aeronautics, vol. 34, no. 4, pp.432-450, 2021, doi: 10.1016/j.cja.2020.07.017.

\bibitem{wall2017survey} T. Wall, R. Meyer, ``A survey of hybrid electric propulsion for aircraft,'' 53rd AIAA/SAE/ASEE Joint Propulsion Conference, pp.4700, 2017, doi: 10.2514/6.2017-4700.

\bibitem{friedrich2015hybrid} C. Friedrich, P.A. Robertson, ``Hybrid-electric propulsion for automotive and aviation applications,'' CEAS Aeronautical Journal, vol. 6, no. 2, pp. 279-290, 2015, doi: 10.1007/S13272-014-0144-X.

\bibitem{brelje2019electric} B.J. Brelje, J.R. Martins, ``Electric, hybrid, and turbo-electric fixed-wing aircraft: A review of concepts, models, and design approaches,'' Progress in Aerospace Sciences, vol. 104, pp. 1-19, 2019, doi: 10.1016/j.paerosci.2018.06.004.

\bibitem{leite2020optimal} J.P.S. Pinto Leite, M. Voskuijl, ``Optimal energy management for hybrid-electric aircraft'', Aircraft Engineering and Aerospace Technology, Vol. 92 No. 6, pp. 851-861, doi: 10.1108/AEAT-03-2019-0046.

\bibitem{doff2020optimal} M. Doff-Sotta, M. Cannon, M. Bacic, ``Optimal energy management for hybrid electric aircraft,'' IFAC-PapersOnLine, vol. 53, no. 2, pp.6043-6049, 2020, doi: 10.1016/j.ifacol.2020.12.1672

\bibitem{donateo2018applying} T. Donateo, A. Ficarella, L. Spedicato, ``Applying dynamic programming algorithms to the energy management of hybrid electric aircraft,'' Turbo Expo: Power for Land, Sea, and Air, vol. 51043, pp. V003T06A015, 2018, doi: 10.1115/GT2018-76500.


\bibitem{martini2011world} F. Martini, ``World’s first serial hybrid electric aircraft to fly at Le Bourget,'' AXX20110666e, Siemens, Munich, 2011.

\bibitem{wittmann2013flying} W. Wittmann, ``Flying with Siemens integrated drive system,'' IDT2013064084e, Siemens, Nuremberg, Germany, 2013.

\bibitem{langelaan2013green} J.W. Langelaan, A. Chakrabarty, A. Deng, K. Miles, V. Plevnik, J. Tomazic, T. Tomazic, G. Veble, ``Green flight challenge: aircraft design and flight planning for extreme fuel efficiency,'' Journal of Aircraft, vol. 50, no.3, pp. 832-846, 2013, doi: 10.2514/1.C032022.


\bibitem{caujolle2017airbus} M. Caujolle, R. Gage, ``Airbus, Rolls-Royce, and Siemens team up for electric future Partnership launches E-Fan X hybrid-electric flight demonstrator,'' 2017. Available: http://www. airbus. com/newsroom/press-releases/en/2017/11/airbus--rolls-royce--and-siemens-team-up-for-electric-future-par. html

\bibitem{de2019preliminary} R. De Vries, M. Brown, R. Vos, ``Preliminary sizing method for hybrid-electric distributed-propulsion aircraft,'' Journal of Aircraft, vol. 56, no. 6, pp. 2172-2188, 2019, doi: 10.2514/1.C035388.

\bibitem{isikveren2014pre} A.R. Isikveren, S. Kaiser, C. Pornet, P.C. Vratny, ``Pre-design strategies and sizing techniques for dual-energy aircraft,'' Aircraft Engineering and Aerospace Technology: An International Journal, 2014.

\bibitem{finger2020comparison} D.F. Finger, R. de Vries, R. Vos, C. Braun, C. Bil, ``A comparison of hybrid-electric aircraft sizing methods,'' AIAA Scitech 2020 Forum, pp. 1006, 2020, doi: 10.2514/6.2020-1006. 


\bibitem{wang2021power} M. Wang, M. Mesbahi, ``Power Allocation For Hybrid Electric Aircraft via Optimal Control during Climb, Cruise, and Descent,'' AIAA Scitech 2021 Forum, pp. 0640, 2021, doi: 10.2514/6.2021-0640.

\bibitem{maurer2005optimization} H. Maurer, I. Altrogge, N. Goris, ``Optimization methods for solving bang-bang control problems with state constraints and the verification of sufficient conditions,'' Proceedings of the 44th IEEE Conference on Decision and Control, pp. 923-928, 2005.

\bibitem{kim2003sensitivity} J.H. Kim, H. Maurer, ``Sensitivity analysis of optimal control problems with bang-bang controls,'' 42nd IEEE International Conference on Decision and Control, vol. 4, pp. 3281-3286, 2003.

\bibitem{maurer1977optimal} H. Maurer, ``On optimal control problems with bounded state variables and control appearing linearly,'' SIAM Journal on Control and Optimization, vol. 15, no. 3, pp. 345-362, 1977.

\bibitem{kaya2003computational} C.Y. Kaya, J.L. Noakes, ``Computational method for time-optimal switching control,'' Journal of optimization theory and applications, vol. 117, no. 1, pp. 69-92, 2003.
\end{thebibliography}
\end{document}